\documentclass[a4paper,12pt]{article}
\def\n{\noindent}
\normalsize
\begin{document}

\title {\bf On greedy algorithms with respect to generalized Walsh system }
\author{S. A. Episkoposian}

\date {e-mail: sergoep@ysu.am}

 \maketitle
{In this paper we proof that there exists a function $f(x)$
belongs to $L^1[0,1]$ such that a greedy algorithm
 with regard to generalized Walsh system does not converge to $f(x)$
in $L^1[0,1]$ norm, i.e. the generalized Walsh system is not a
quasi-greedy basis in its linear span $L^1[0,1]$.} \vskip 4mm

{\it Keywords}: greedy algorithms, generalized Walsh system,
quasi-greedy basis.

{\it 2000 Mathematics Subject Classification}: Primary  42C10;
Secondary 46E30.

\section{\bf Introduction }
\par\par\bigskip
\par\par\bigskip

 In this paper we consider a question of convergence of
greedy algorithm with regard to generalized Walsh system in
$L^1[0,1]$ norm.

Let $X$ be a Banach space  with a norm $||\cdot||=||\cdot||_X$ and
a basis $\Phi=\{\phi_k\}_{k=1}^\infty$, $||\phi_k||_X=1$,
$k=1,2,..$ .

Denote by $\Sigma_m$ the collection of all functions in $X$ which
can be expressed as a linear combination of at most $m$- functions
of $\Phi$. Thus each function $g\in \Sigma_m$ can be written in
the form
$$g=\sum_{s\in \Lambda}a_s\phi_s,\ \ \#\Lambda \leq m.$$

For a function $f\in X$ we define its approximate error by
$$\sigma_m(f,\phi)=\inf_{g\in \Sigma_m}\left| \left| f-g \right| \right|_X,\ \ m=1,2,...$$

and we consider the expansion
$$f=\sum_{k=1}^\infty a_k(f)\phi_k \ \ .$$

{\bf Definition 1.} Let an element $f\in X$ be given. Then the
$m$-th greedy approximant of the function $f$ with regard to the
basis $\Phi$ is given by
$$G_m(f,\phi)=\sum_{k\in\Lambda}a_k(f) \phi_k,\eqno(1)$$
where $\Lambda \subset \{1,2,...\}$ is a set of cardinality $m$
such that
 $$|a_n(f)|\geq |a_k(f)|, \ \  n\in \Lambda,\ \  k\notin
\Lambda, \eqno (2)$$ We'll say that the greedy approximant of
$f(t)\in L^p_{[0,1]}$, $p\geq 0$ converges with regard to the
basis $\Phi$ , if the sequence $G_m(x,f)$ converges to $f(t)$ in
$L^p$ norm. This new and very important direction invaded many
mathematician's attention (see [1]-[10]).

{\bf Definition 2.} We call a basis $\Phi$ greedy basis if for
every $f \in X$  there exists a subset $\Lambda\subset
\{1,2,...\}$ of cardinality $m$, such that
$$\left| \left| f-G_m(f,\Phi) \right| \right|_X\leq C\cdot\sigma_m(f,\Phi)$$
where a constant $C=C(X,\Phi)$ independent of $f$ and $m$.

In [4] it is proved that each basis $\Phi$ which is $L_p$
-equivalent to the Haar basis $H$ is Greedy basis for $L_p(0,1)$,
$1<p<\infty$ .

{\bf Definition 3.} We say that a basis $\Phi$ is Quasi-Greedy
basis if there exists a constant $C$  such that for every $f\in X$
and any finite set of indices $\Lambda$, having the property
$$\min_{k\in \Lambda}|a_k(f)|\geq \max_{n\notin \Lambda}|a_k(f)|$$
we have
$$||S_{\Lambda}(f,\Phi)||_X=\left| \left| \sum_{k\in \Lambda}a_k(f)\phi_k \right| \right|_X \leq C\cdot ||f||_X.$$
\vskip 4mm

In [5] it is proved that a basis $\Psi$ is quasi-greedy if and
only if the sequence $\{G_m(f) \}$ converges to $f$, for all $f\in
X$. Not that the trigonometric and Walsh system are not a
quasi-greedy basis for $L^p$ if $1<p<\infty$ (see [7] and [8] ).
\par\par\bigskip
\par\par\bigskip

\section{\bf  Definition and properties of generalized Walsh system }
\par\par\bigskip

Let $a$ denote a fixed integer, $a\geq 2$ and put $\omega_a=e^{2
\pi i \over a}$.

Now we will give the definitions of Rademacher and generalized
Walsh systems (see [12] ).

{\bf Definition 4.} The Rademacher  system of order $a$ is defined
by
$$\varphi_0(x)=\omega_a^k\ \ if \ \ x \in \left[ {k\over a}, {k+1\over a}\right),\ \ k=0,1,...,a-1,$$
and for $n\geq 0$
$$\varphi_n(x+1)=\varphi_n(x)=\varphi_0(a^nx).\eqno(3)$$

{\bf Definition 5.} The generalized Walsh system of order $a$ is
defined by
$$\psi_0(x)=1,$$
and if $n=\alpha_{n_1}a^{n_1}+...+\alpha_{n_s} a^{n_s}$ where
$n_1>...>n_s,$ then
$$\psi_n(x)=\varphi_{n_1}^{\alpha_{n_1}}(x)\cdot...\cdot \varphi_{n_s}^{\alpha_{n_s}}(x) .\eqno(4)$$
Let $\Psi_{a}=\{ \psi_n(x)\}_{n=0}^\infty$ denote the generalized
Walsh system of order $a$. Not that $\Psi_2$ is the classical
Walsh system.

{\bf Remark.} The generalized Walsh system $\Psi_a$, $a\geq 2$ is
a complete orthonormal system in $L^2[0,1)$ (see [12]).

The basic properties of the generalized Walsh system of order $a$
are obtained by  R. Paley , H.E.Chrestenson, J. Fine, N. Vilenkin
and others (see [11]- [16]).

 Define
$$I_{n,k}=I_{n,k}(a)= \left[ {k\over a^n}, {k+1\over a^n}\right),\ \ k=0,...,a^n-1,\ \ n=1,2,....$$
If $\varphi_n(x)$ is the $n$th Rademacher function of order $a$,
then from Definition 4 it follows
$$\varphi_n(x)=\omega_a^k=e^{2 \pi
i \cdot k \over a},\ \ x\in I_{n+1,k}.\eqno(5)$$ Note some
properties of generalized Walsh system: \vskip 4mm

{\it Property 1.} From definition 5 we have
$$\psi_{a^k+j}(x)=\varphi_k(x) \cdot \psi_j(x),\ \ if\ \ 0 \leq j \leq a^k-1.\eqno (6)$$
Denote by
$$D_n(t)=\sum_{k=0}^{n-1}\psi_k(t), \eqno(7)$$
the Dirichlet kernel by generalized Walsh system.

{\it Property 2.} The Dirichlet kernel has the following
properties (see [12] )

$$D_{a^n}(t)=\cases{a^n , \ \ x \in I_{n,0}=\big[0,{1\over a^n}\big);\cr \cr 0 ,
\ \ x \in \big[{1\over a^n},1\big).\cr}\eqno (8)$$

{\it Property 3.} If $n=a^k+m$, $0\leq m<k$ and consequently by
(6) - (8) we have
$$D_n(t)=D_{a^k}(t)+\varphi_k(t)\cdot D_m(t),\ \ t\in [0,1). \eqno(9)$$

{\it Property 4.} For any natural number $m$ and any $t\in (0,1)$
the following is true
$$|D_m(t)|\leq m . \eqno(10)$$
\par\par\bigskip
\par\par\bigskip

\section {\bf A Basic Lemma }
\par\par\bigskip

Denote by
$$L_k=\int_0^1|D_k(t)|dt$$
the $k$th Lebesgue constant of the generalized Walsh system
$\{\Psi_a\}$.

 In [12] it is proved that the Lebesgue constant satisfy $L_k=O(\log_a k)$
where $O$ depends upon $a$. Next Lemma shows that there exists a
sequence of natural numbers $\{n_k\}$ so that the sequence
$L_{n_k}$ has the same order of growth as $\log_a n_k$. Namely the
following is true: \vskip 4mm

{\bf Lemma .} There exists a sequence of natural numbers
$\{n_k\}_{k=1}^\infty$ of the form
$$n_{2s}=\sum_{i=0}^s a^{2i}:\ \ n_{2s+1}=\sum_{i=0}^s a^{2i+1},\ \ s=0,1,2,...,\eqno(11)$$
such that
$$a^k\leq n_k<a^{k+1},$$
$$L_{n_k}=\int_0^1|D_{n_k}(t)|dt>{1 \over a}\cdot \left({k\over 2}+
1 \right)>{1\over 2a}\cdot \log_a n_k,\ \ k\geq 1. \eqno(12)$$
\par\par\bigskip

{\bf Proof.} Note that

$$n_{2s}={a^{2s+2}-1 \over a^2-1}<{a^2 \over a^2-1}\cdot a^{2s},$$
$$n_{2s+1}={a^{2s+3}-a \over a^2-1}=a\cdot {a^{2s+2}-1 \over a^2-1} <{a^2 \over a^2-1}\cdot a^{2s+1}$$
i.e.
$$n_k<{a^2 \over a^2-1}\cdot a^k<a^{k+1},\ \ k=0,1,2,... \eqno(13)$$
From this and (10) we have
$$|D_{n_k}(t)|<{a^2 \over a^2-1}\cdot a^k,\ \ t\in [0,1), \ \ k=0,1,2,... \eqno(14)$$
First we'll prove that
$$\int_{1 \over a^{k+2}}^1|D_{n_k}(t)|dt>{1 \over a}\cdot \left({k\over 2}+
1 \right)\eqno(15)$$

Let $k=2s$, then we have to prove
$$\int_{1 \over a^{2s+2}}^1|D_{n_{2s}}(t)|dt>{1 \over a}\cdot \left(s+1 \right),\ \ s=0,1,2,...\eqno(16)$$
By Definition 5 and (7) for $s=0$ we have
$$\int_{1 \over a^2}^1|D_1(t)|dt=\int_{1 \over a^2}^1|\psi_0(t)|dt=1-{1 \over a^2}>{1 \over
a}\ \ .$$
 Now assume that for some $s-1$ the inequality (16)
holds, i.e.
$$\int_{1 \over a^{2s}}^1|D_{n_{2(s-1)}}(t)|dt>{1 \over a}\cdot s .\eqno(17)$$
By (8) and (14) we get
$$D_{a^{2s}}(t)=a^{2s},\ \ if \ \ t\in I_{2s,0}.\eqno(18)$$
$$|D_{n_{2(s-1)}}(t)|<{a^2 \over a^2-1}\cdot a^{2(s-1)}={a^{2s} \over a^2-1}.\eqno(19)$$
From (11) it follows that
$$n_{2s}=a^{2s}+n_{2(s-1)}\eqno(20)$$
and consequently by (9), (17) and (18) for $t\in I_{2s,0}$ we have
$$|D_{n_{2s}}(t)|=|D_{a^{2s}}(t)+\varphi_{2s}(t)\cdot D_{n_{2(s-1)}}(t)|\geq$$
$$|D_{a^{2s}}(t)|-|D_{n_{2(s-1)}}(t)|> {a^2-2 \over a^2-1}\cdot a^{2s}.$$
Hence, taking into account that $\left({1 \over a^{2s+2}},{1 \over
a^{2s}}\right)\subset I_{2s,0}=\left[0,{1 \over a^{2s}}\right)$,
we obtain
$$\int_{1 \over a^{2s+2}}^{1 \over a^{2s}}|D_{n_{2s}}(t)|dt>{a^2-2 \over a^2-1}\cdot a^{2s}\cdot
\left({1 \over a^{2s}}-{1 \over a^{2s+2}}\right)=$$
$${a^2-2 \over a^2-1}\cdot {a^2-1\over a^{2s+2}}={a^2-2 \over a^2}\geq {1 \over a}.\eqno(21)$$
From (8), (9), (17) and (20) follows
$$\int_{1 \over a^{2s}}^1|D_{n_{2s}}(t)|dt>
\int_{1 \over a^{2s}}^1|D_{n_{2(s-1)}}(t)|dt>{1 \over a}\cdot s.
$$
Hence and from (21) we conclude
$$\int_{1 \over a^{2s+2}}^1|D_{n_{2s}}(t)|dt=$$
$$\int_{1 \over a^{2s+2}}^{1 \over a^{2s}}|D_{n_{2s}}(t)|dt+
\int_{1 \over a^{2s}}^1|D_{n_{2s}}(t)|dt>$$
$${1 \over a}+{1 \over a}\cdot s> {1 \over a}\cdot \left(s+1
\right),\ \ s=0,1,2,...\eqno(22)$$

In a case $k=2s+1$ $(s=0,1,...)$ we have to prove
$$\int_{1 \over a^{2s+3}}^1|D_{n_{2s+1}}(t)|dt>{1 \over a}\cdot \left(s+{1\over 2}+1 \right).$$
For $s=0$ this inequality holds because in this case
$n_{2s+1}=n_1=a$ and
$$\int_{1 \over a^3}^1|D_a(t)|dt=\int_{1 \over a^3}^{1 \over a}adt=a\cdot \left({1 \over a}-{1 \over
a^3}\right)=$$
$${1 \over a} \cdot \left(a-{1 \over a}\right)\geq
{1\over a}\cdot {3\over 2}.$$
The next  reasonings are similar to
a case when $k=2s$.

Since $a^k\leq n_k<a^{k+1}$ then $\log_an_k<k+1$ and consequently
$$L_{n_k}=\int_0^1|D_{n_k}(t)|dt>\int_{1 \over a^{k+2}}^1|D_{n_k}(t)|dt>{1 \over a}\cdot \left({k\over 2}+
1 \right)>$$
$${1\over 2a}\cdot \left({k \over 2}+1\right)>{1\over 2a}\cdot\log_a n_k.$$

\par\par\bigskip
{\bf Completing the proof.}
\par\par\bigskip
\par\par\bigskip
\break

\section{\bf The Main Theorem and It's Proof.}
\par\par\bigskip

In [10] we proved the following theorem: \vskip 4mm
 {\bf Theorem 1.} Let a sequence $\{ M_n \}_{n=1}^\infty$ be given so that
$$\lim_{k \to \infty}(M_{2k}-M_{2k-1})=+\infty.$$
Then the Walsh subsystem
$$\{ W_{n_k}(x) \}_{k=1}^\infty=\{W_ m(x):\ \ M_{2s-1}\leq m \leq M_{2s},\ \ s=1,2,...\}\eqno(1)$$
is not a quasi-greedy basis in its linear span in $L^1[0,1]$.
\vskip 4mm
 Let $\Psi_{a}=\{\psi_n(x)\}_{n=0}^\infty$ denote the generalized Walsh system of
order $a$.

 From Corollary 2.3 (see [8]) it follows that generalized Walsh
 system is not a quasi-greedy basis for $L^p[0,1]$ if
 $1<p<\infty$.

In this paper we prove the following theorem. \vskip 4mm

{\bf Theorem 2.}  There exists a function $f(x)$ belongs to
$L^1[0,1)$ such that the approximate $G_n(f,\Psi_a)$ with regard
to the generalized Walsh system does not converge to $f(x)$ by
$L^1[0,1)$ norm, i.e. the generalized Walsh system is not a
quasi-greedy basis in its linear span in $L^1$.\vskip 4mm

{\bf Proof.} Let $a\geq 2$ denote a fixed integer. For any natural
$k$ we set

$$f_k(x)=\sum_{i=a^{(k-1)^2}}^{a^{k^2}-1}\left( {1 \over k^2}+2^{-i}\right) \cdot \psi_i(x).\eqno(23)$$
It is easy to see that the Fourier coefficients by generalized
Walsh system of the function $f_k(x)$ are defined as follows
$$C_i^{(k)}={1 \over k^2}+2^{-i}\ \ if \ \ a^{(k-1)^2}\leq i<a^{k^2}.\eqno(24)$$
Now we consider the following function
$$f(x)=\sum_{i=1}^\infty C_i \psi_i(x)=\sum_{k=1}^\infty f_k(x)=$$
$$=\sum_{k=1}^\infty\left[ \sum_{i=a^{(k-1)^2}}^{a^{k^2}-1}\left( {1 \over
k^2}+2^{-i}\right) \cdot \psi_i(x)\right ],\eqno(25)$$
where
$$C_i=C_i^{(k)}\ \ if \ \ a^{(k-1)^2}\leq i<a^{k^2}.\eqno(26)$$

Now we will show that $f(x)\in L^1[0,1]$. For this we represent
the function $f(x)$ in the following way:
$$f(x)=g(x)+h(x),\eqno(27)$$
where
$$g(x)=\sum_{k=1}^\infty {1 \over
k^2} \left[ \sum_{i=a^{(k-1)^2}}^{a^{k^2}-1} \psi_i(x)\right
]=\sum_{k=1}^\infty {1 \over k^2} \left[
D_{a^{k^2}}(x)-D_{a^{(k-1)^2}}(x)\right ],$$

$$h(x)=\sum_{k=1}^\infty \left[ \sum_{i=a^{(k-1)^2}}^{a^{k^2}-1} 2^{-i}\cdot \psi_i(x)\right ]=
\sum_{j=1}^\infty 2^{-j}\cdot \psi_j(x).$$

For the function $g(x)$ and from (8) and definition 5 we have
$$\int_0^1|g(x)|dx\leq 2\cdot \sum_{k=1}^\infty {1 \over k^2}<\infty$$
which means $g(x)\in L^1[0,1]$.

Analogously
$$\int_0^1|h(x)|dx\leq \sum_{j=1}^\infty {1 \over 2^j}<\infty$$
i.e. $h(x)\in L^1[0,1]$. Hence and from (27) it follows that
$f(x)\in L^1[0,1]$.

For any natural $k$ we choose numbers $i,j$ so that
$$a^{(k-1)^2}\leq i <a^{k^2}\leq j <a^{(k+1)^2}.$$
Then
$${1 \over (k+1)^2}+2^{-j}<{1 \over k^2}+2^{-i},$$
and from (24) we have $C_j^{(k+1)}(f)<C_i^{(k)}(f)$.

Analogously for any number $i$, $a^{(k-1)^2}\leq i <a^{k^2}$, we
have
$${1 \over k^2}+2^{-(i+1)}<{1 \over k^2}+2^{-i},$$
i.e. $C_{i+1}^{(k)}(f)<C_i^{(k)}(f).$

Thus for any natural numbers $i$ we get
$$C_{i+1}(f)<C_i(f).$$
On the other hand if $i\to \infty$ then $k\to \infty$ (see (24)).
Then from (24) and (26) we get $C_i(f) \searrow 0$.

For any numbers $m_k $ so that
$$a^{(k-1)^2}+ m_k<a^{k^2},\eqno(28)$$
by (24) - (26) and Definition 1 we have
$$G_{a^{(k-1)^2}+ m_k}(f,\Psi_a)-G_{a^{k^2}}(f,\Psi_a)= $$
$$\sum_{i=a^{(k-1)^2}}^{a^{(k-1)^2}+m_k-1}C_i^{(k)} \cdot \psi_i(x)=
{1 \over k^2} \sum_{i=a^{(k-1)^2}}^{a^{(k-1)^2}+m_k-1}\psi_i(x)+$$
$$ \sum_{i=a^{(k-1)^2}}^{a^{(k-1)^2}+m_k-1} {1 \over 2^i}\cdot \psi_i(x)=J_1+J_2. \eqno(29)$$

Taking into account (6) and (7) we get
$$J_1={1 \over k^2} \sum_{i=a^{(k-1)^2}}^{a^{(k-1)^2}+m_k-1}\psi_i(x)
={1 \over k^2} \sum_{i=0}^{m_k-1}\psi_{a^{(k-1)^2}+i}(x)=$$
$${1 \over k^2}\cdot \psi_{a^{(k-1)^2}}(x)\cdot \sum_{i=0}^{m_k-1}\psi_i(x)
={1 \over k^2}\cdot \psi_{a^{(k-1)^2}}(x)\cdot D_{m_k}(x). $$

$$\mid J_2 \mid \leq \sum_{i=a^{(k-1)^2}}^{a^{(k-1)^2}+m_k-1} {1 \over 2^i} \mid \psi_i(x)\mid
\leq \sum_{i=a^{(k-1)^2}}^{\infty} {1 \over 2^i}\leq
2^{-a^{(k-1)^2}+1}.$$

From this and (29) we obtain\vskip 4mm
$$\mid G_{a^{(k-1)^2}+ m_k}(f,\Psi_a)-G_{a^{k^2}}(f,\Psi_a) \mid \geq $$
$$ {1 \over k^2}\cdot \mid \psi_{a^{(k-1)^2}}(x) \mid \mid D_{m_k}(x)\mid-2^{-a^{(k-1)^2}+1}=$$
$${1 \over k^2}\cdot \mid D_{m_k}(x)\mid-2^{-a^{(k-1)^2}+1}.\eqno(30)$$

Now we take the sequence of natural numbers $m_\nu$ defined by
Lemma  (see (11) è (12))such that $a^{(k-1)^2}\leq m_\nu <
a^{(k-1)^2+1}$.

 Then from (30) we have\vskip 4mm
$$\int_0^1\mid G_{a^{(k-1)^2}+ m_k}(f,\Psi_a)-G_{a^{k^2}}(f,\Psi_a)\mid dx>$$
$${1 \over k^2}\cdot \int_0^1|D_{m_k}(x)|dx - 2^{-a^{(k-1)^2}+1}
\geq{1\over 2a\cdot k^2}\cdot \log_a m_k - 2^{-a^{(k-1)^2}+1} \geq $$
$$ {(k-1)^2 \over 2a \cdot k^2} - 2^{-a^{(k-1)^2}+1}
\geq {1 \over 4a}- 2^{-a^{(k-1)^2}+1} \geq C_1,\ \ \ \ k \geq
4.$$\vskip 4mm

Thus  the sequence $\{G_n(f,\Psi_a)\}$ does not converge  by
$L^1[0,1]$ norm, i.e. the generalized Walsh system $\Psi_a$ is not
a quasi-greedy basis in its linear span in $L^1[0,1]$.

\par\par\bigskip
{\bf Completing the proof.}
\par\par\bigskip
\par\par\bigskip
\break

 {\bf R E F E R E N C E S }

\par\par\bigskip
\par\par\bigskip
\par\par\bigskip

\n [1] DeVore R. A., Temlyakov V. N., Some remarks on Greedy
Algorithms, Adv. Comput. Math., 1995, v.5, p.173-187.
\par\par\bigskip

\n [2] DeVore  R. A., Some remarks on greedy algorithms. Adv
Comput. Math., 1996, v.5, p. 173-187.
\par\par\bigskip

\n [3] Davis G., Mallat S. and Avalaneda M., Adaptive greedy
approximations. Constr, Approx., 1997, v.13, p. 57-98.
\par\par\bigskip

\n [4] Temlyakov V. N., The best $m$ - term approximation and
Greedy Algorithms,  Advances in Comput. Math., 1998, v.8,
p.249-265.
\par\par\bigskip

\n [5] Wojtaszcyk P., Greedy Algorithm for General Biorthogonal
Systems, Journal of Approximation Theory, 2000, v.107, p. 293-314.
\par\par\bigskip

\n [6] Konyagin S. V., Temlyakov V. N., A remark on Greedy
Approximation in Banach spaces, East Journal on Approximation,
1999, v.5, p. 493-499.
\par\par\bigskip

\n [7] Temlyakov V. N., Greedy Algorithm and $m$ - term
trigonometric approximation, Constructive Approx., 1998, v.14, p.
569-587.
\par\par\bigskip

\n [8] Gribonval R., Nielsen M., On the quasi-greedy property and
uniformly bounded orthonormal systems ,
http://www.math.auc.dk/research/reports/R-2003-09.pdf.
\par\par\bigskip

\n [9] Grigorian M.G., "On the convergence of Greedy algorithm",
International Conference, Mathematics in Armenia, Advances and
Perspectives,Yerevan, Armenia, Abstract, 2003,  p. 44 - 45.
\par\par\bigskip

\n [10] Episkoposian S.A., On the divergence of Greedy algorithms
with respect to Walsh  subsystems in $L^1$ ,  Nonlinear Analysis:
Theory, Methods \& Applications, 2007, v. 66, Issue 8, p.
1782-1787.
\par\par\bigskip

\n [11] Pely R.,  A remarkable systems of orthogonal functions,
Proc. London Math. Soc., 1932, v.34, p. 241-279.
\par\par\bigskip

\n [12] Chrestenson H.E., A class of generalized Walsh functions,
Pacific J. Math., 1955, v.45, p. 17-31.
\par\par\bigskip

\n [13] Fine J., The generalized Walsh functions Trans. Amer.
Math. Soc., 1950, v.69, p. 66-77.
\par\par\bigskip

\n [14] Vilenkin N., On a class of complete orthonormal systems,
AMS Transl., 1963, v.28, p. 1-35.
\par\par\bigskip

\n [15] Selfridge R., Generalized Walsh transform, Pacific J.
Math., 1955, v.5, p. 451-480.
\par\par\bigskip

\n [16] Watari C., On generalized Walsh-Fourier series, Tohoku
Math. J., 1958, v.10, p. 211-241.
\par\par\bigskip

\par\par\bigskip
\par\par\bigskip
\par\par\bigskip

Department of Mathematics,

State Engineering University of Armenia

Yerevan, Terian st. 105,

e-mail: sergoep@ysu.am

\end{document}